\documentstyle[12pt]{article}
\begin{document}
\newtheorem{thm}{Theorem}
\newtheorem{prop}{Proposition}
\newtheorem{cor}{Corollary}
\newtheorem{ex}{Example}
\newtheorem{exs}{Examples}
\newtheorem{lem}{Lemma}
\newtheorem{rem}{Remark}
\newcommand{\bt}{\begin{thm}}
\newcommand{\et}{\end{thm}}
\newcommand{\bl}{\begin{lem}}
\newcommand{\el}{\end{lem}}
\newcommand{\bp}{\begin{prop}}
\newcommand{\ep}{\end{prop}}
\newcommand{\bc}{\begin{cor}}
\newcommand{\ec}{\end{cor}}
\newcommand{\p}{{{\bf Proof.\,\,}}}
\def\NN{I\!\! N}
\def\RR{I\!\! R}
\def\HH{I\!\! H}
\def\QQ{I\!\!\!\! Q}
\def\CC{I\!\!\!\! C}
\def\ZZ{Z\!\!\! Z}

\newcommand{\mod}[3]{\mbox{${#1} \equiv{#2}\bmod{#3}$}}
\newcommand{\G} {{\bf G}}
\newcommand{\biindice}[3]%
{
\renewcommand{\arraystretch}{0.5}
\begin{array}[t]{c}
#1\\
{\scriptstyle #2}\\
{\scriptstyle #3}
\end{array}
\renewcommand{\arraystretch}{1}
}
\newcommand{\chit}{\mbox{${\tilde\chi}$}}
\def\qed{\hfill \vrule height4pt width4pt depth2pt}
\baselineskip6mm

\centerline{{\Large Centralizing automorphisms and Jordan left
derivations  }}
 \centerline{{\Large on $\sigma$-prime rings }} \vspace*{0,6cm}
\centerline{{\large L. Oukhtite \& S. Salhi }}
\centerline{\small{Universit\'e Moulay Isma\"\i l, Facult\'e des
Sciences et Techniques}} \centerline{\small{D\'epartement de
Math\'ematiques, Groupe d'Alg\`ebre et Applications}}
\centerline{\small{ B. P. 509 Boutalamine, Errachidia; Maroc}}
\centerline { oukhtitel@hotmail.com, \,\,\,salhi@math.net}
\begin{abstract}\vspace*{-0,5cm}
\noindent Let $\;R\;$ be a $2$-torsion free $\sigma$-prime ring. It
is shown here that if  $U\not\subset Z(R)$ is a $\sigma$-Lie ideal
of $R$ and $a, b$ in $R$ such that $aUb=\sigma(a)Ub=0,$ then either
$a=0$ or $b=0.$ This result is then applied to
 study the relationship between the structure of
R and certain automorphisms on R. To end this paper, we describe
additive maps $d: R \longrightarrow R$ such that $d(u^2) = 2ud(u)$
where $u\in U,$ a nonzero $\sigma$-square closed Lie ideal of $R.$
\end{abstract}
\vspace*{-0,4cm}
\section{{\small{\bf Introduction}} }
\vspace*{-0,4cm} \indent Recently, there has been much interest in
investigating the structure of a ring having commuting or
centralizing mapping on a subset $S$ of this  ring. In \cite{Pos}
Posner showed that if a prime ring has a nontrivial derivation which
is centralizing on the entire ring, then the ring must be
commutative. In \cite{May76} the same result is proved for a prime
ring with a nontrivial centralizing automorphism. A number of
authors have generalized these results by considering mappings which
are only assumed to be centralizing on an appropriate ideal of the
ring. In \cite{May84} Mayne showed that if $R$ is a prime ring with
a nontrivial centralizing automorphism on a nonzero ideal or
(quadratic) Jordan ideal, then R is commutative. For a prime ring
with characteristic not two, an improved version was given in
\cite{May92} by showing that a nontrivial automorphism which is
centralizing on a Lie ideal implies that the ideal is contained in
the center of the ring. In \cite{O.S2}, we proved a corresponding
result for $\sigma$-prime rings with $\sigma$-square closed Lie
ideals, where $\sigma$ is an involution. More precisely, for  a
$\sigma$-prime ring $R$ with characteristic not $2$ equipped with a
nontrivial automorphism $T$ centralizing on a $\sigma$-square closed
Lie ideal $U$, it is proved that if $T$ commutes with $\sigma$ and
there exists $u_0$ in $Sa_{\sigma}(R)$ with $Ru_0\subset U,$ then
$U$ is contained in $Z(R).$ This extends the result of \cite{May84}
for centralizing automorphisms to $\sigma$-prime rings.\\Our first
aim in this paper is to give a fundamental lemma which shows that if
$U$ is a noncentral $\sigma$-Lie ideal of $R$ and $a, b$ in $R$ such
that $aUb=\sigma(a)Ub=0,$ then $a=0$ or $b=0$, Lemma 4. The
application of this lemma enables us to improve the result of
\cite{O.S2}  by showing that the existence of $u_0$ in
$Sa_{\sigma}(R)$ such that $Ru_0\subset U$ is not necessary, Theorem 1.\\
 \indent To end this paper, we present another field of application of Lemma 4
 by  studying the relationship between Jordan
 left derivations and left derivations on
 Lie ideals. A famous result due to Herstein \cite{Her}  states that every Jordan
derivation on a 2-torsion free prime ring is a derivation. Further,
Awtar \cite{Awt} generalized this result on square closed Lie
ideals. In \cite{Ash} Mohammad Ashraf and Nadeem -UR-Rehmman proved
that if $d$ is an additive mapping of a $2$-torsion free prime ring
$R$ into itself satisfying $d(u^2)=2ud(u),$ for all $u$ in a square
closed Lie ideal $U$ of $R,$ then $d(uv)=ud(v)+vd(u)$ for all
$u,v\in U.$ Our main goal in Theorem 2 will be to  extend this
result to $\sigma$-prime rings. \vspace*{-0,4cm}
\section{ {\small {\bf Preliminaries and notations}}}
\vspace*{-0,4cm}In the starting of this paper we should give some
definitions and preliminary results to be used. Throughout, $R$ will
represent an associative ring with center $Z(R).$ $R$ is said to be
$2$-torsion free if whenever $2x=0,$ with $x\in R,$ then $x=0.$ As
usual the commutator $xy-yx$ will be denoted by $[x,y].$ We shall
use basic commutator identities $\;[x, yz] = y[x, z] + [x, y]z,\;\;
[xy,z] = x[y,z] + [x,z]y .$ An involution $\sigma$ of a ring $R$ is
an anti-automorphism of order $2$ (i.e. $\sigma$ is an additive
mapping satisfying $\sigma(xy)= \sigma(y)\sigma(x)$ and $\sigma
^2(x)=x$ for all $x,y\in R).$ Given an involutorial ring  $(R,
\sigma),$  we define $\;Sa_{\sigma}(R) :=\{ r\in \;\;R /
\;\sigma(r)=\pm \, r \}.$\\Recall that $R$ is $\sigma$-prime if $aRb
= aR\sigma(b)=0$ implies that $a = 0$ or $b = 0.$ It is worthwhile
to note that every prime ring having an involution $\sigma$ is a
$\sigma$-prime ring but the converse is in general not true
(\cite{Diaspo}). A Lie ideal of $R$ is an additive subgroup $U$ of
$R$ satisfying $[U,R]\subset U$. A Lie ideal $U$ of $R$ is said to
be a $\sigma$-Lie ideal, if $\sigma(U)=U.$ If $U$ is a $\sigma$-Lie
ideal of $R$ such that $u^2 \in U$ for all $u\in U,$ then $U$ is
called a $\sigma$-square closed Lie ideal. For $u, v \in U, (uv +
vu) = (u + v)^2 - (u^2 + v^2)$ and so $uv + vu \in U.$ Also, $[u, v]
= uv - vu \in U$ and it follows that $2uv \in U.$ This remark will
be freely used in the whole paper. A linear mapping $T$ of a ring
$R$ into itself is called centralizing on a subset $S$ of $R$ if
$[x,T(x)]\in Z(R)$ for every $x$ in $S.$ In particular, if $T$
satisfies $[x, T(x)] = 0$ for all
$x$ in $S$ then $T$ is called commuting on $S.$\\
An additive mapping $d : R\longmapsto R$ is called a left derivation
(resp. Jordan left derivation)   if $d(xy)= xd(y) + y d(x)$ (resp.
$d(x^2)=2xd(x)$) holds for all pairs $x,y\in R.$ Clearly, every left
derivation is a Jordan left derivation, but the converse need not be
true in general. \vspace*{-0,4cm}
\section{ {\small {\bf  Centralizing automorphisms on $\sigma$-square closed Lie ideals}} }
Let $R$ be a prime ring of characteristic not $ 2$ and $U$ a Lie
ideal of $R$ with $U \not\subset Z(R).$ Bergen et al. proved in
(\cite{Ber}, Lemma 4) that if $a,b \in R$ are such that $aUb = 0$,
then either $a = 0$ or $b = 0.$ This Lemma is the key of the
intensive study of the relationship  between several maps
(especially derivations and automorphisms) and Lie ideals of prime
rings. In particular, when this relationship involves the action of
these maps on Lie ideals. Many of these results extend other ones
proven previously just for the action of the maps on the whole ring.
In this section, we are primarily interested in centralizing
automorphisms on Lie ideals so that our first object will be the
extension of Lemma 4 of \cite{Ber} to $\sigma$-prime rings. This
extension will also play an important
role in the last section of the present paper.\\
In what follows $R$ will denote a $\sigma$-prime ring of
characteristic not $2.$
\vspace*{-0,3cm}
\bl\label{l1}

{\em (\cite{O.S2}, Lemma 4)} Let $0\neq I$ be a $\sigma$-ideal of
a $2$-torsion free $\sigma$-prime ring $R.$ If $a, b$ in $R$ are
such that $aIb=0=aI\sigma(b),$ then $a=0$ or $b=0.$

\el
\vspace*{-0,5cm}
\bl\label{l3}

Let $0\neq U$ be a $\sigma$-Lie ideal of a $2$-torsion free
$\sigma$-prime ring $R.$ If $[U,U]= 0,$ then $U \subset Z(R).$

\el \vspace*{-0,5cm} \p Let $u\in U\cap Sa_\sigma(R),$  from
$[U,U]=0$ it follows that $[u,[u,x]]=0$  for all $x\in R.$ Hence
$d_u^2(x)=0$ for all $x\in R,$ where $d_u$ is the inner derivation
defined by $d_u(y)=[u,y].$ Replacing $x$ by $xy$, we get
$d_u(x)d_u(y)=0$ for all $x,y\in R$ because $char R\neq 2$. If we
replace $x$ by $xz$ in the last equality we obtain
$d_u(x)zd_u(y)=0$ for all $x,y,z\in R$, so that $d_u(x)Rd_u(y)=0$.
As $d_u\circ \sigma=\pm \sigma\circ d_u$ then
$\;d_u(x)Rd_u(y)=0=\sigma(d_u(x))Rd_u(y)\;$ which proves $d_u=0$
in such a way that $u\in Z(R)$. Consequently, $U\cap
Sa_\sigma(R)\subset Z(R).$
   Let $x\in U,$ as $x + \sigma(x)$ and $x-\sigma(x)$ are in  $U\cap Sa_\sigma(R),$
   then $x + \sigma(x)$ and $x-\sigma(x)$ are in $Z(R)$ so that
   $2x\in Z(R).$ Accordingly, $x$ in $Z(R)$ proving $U\subset
   Z(R).$\qed
   \vspace*{-0,3cm}
 \bl\label{l4}

Let $0\neq U$ be a $\sigma$-Lie ideal of a $2$-torsion free
$\sigma$-prime ring $R.$ If $[U,U]\neq 0$, then there exists a
$\sigma$-ideal $M$ of $R$ such that $[M,R]\subset U$ with
$[M,R]\not\subset Z(R)$.

\el \p We have $\;R[U,U]R\subset T(U)\;$ by [2, Lemma 1]  \,where \\
$T(U)=\{x\in R\,:\, [x,R]\subset U\}$. Set $M=R[U,U]R,$ clearly $M$
is a nonzero $\sigma$-ideal of $R$ such that $[M,R]\subset U$ and
$[M,R]\not\subset Z(R).$ For if $[M,R]\subset Z(R),$ then
$0=[[m,rm],r]=[[m,r]m, r]=[m, r]^2.$ As $[m, r]\in Z(R),$ this
yields $[m, r]=0$ so that $m\in Z(R)$ and thus $M\subset Z(R).$ For
any $r,s\in R,$ we have
\vspace*{-0,3cm}$$rsm=mrs=(mr)s=s(mr)=srm\;\;\mbox{ for all}\;\;m\in
M.\vspace*{-0,3cm}$$ Hence $[r,s]m=0$ thereby $[r,s]M=0.$ Applying
Lemma \ref{l1}, we find that $[r,s]=0$ and
$R$ is then commutative which contradicts $[U, U]\neq 0.$\qed\\
\\
All is prepared to introduce our main lemma which extends Lemma 4 of
\cite{Ber} to $\sigma$-prime rings.
\vspace*{-0,3cm}
 \bl\label{l5}

If $U\not\subset Z(R)$ is a $\sigma$-Lie ideal of a $2$-torsion
free $\sigma$-prime ring $R$ and $a,b\in R$ such that
$aUb=\sigma(a)Ub=0,$ then $a=0$ or $b=0.$

\el
\vspace*{-0,5cm}
\p Suppose $b\neq 0$, we then conclude that
$Ub\neq 0.$ Otherwise, for any $x\in U$ we get $[x,r]b=0$ for all
$r\in R.$ Replace $r$ by $rt,$ where $t\in R,$ to get $0=[x,
rt]b=[x,r]tb$ so that $[x, r]Rb=0$ for all $r\in R.$ As
$\sigma(U)=U,$ it follows that $[x, r]Rb =0 = \sigma([x, r])R b$ and
the $\sigma$-primeness of $R$ yields $[r, x]=0$ so that $U \subset
Z(R)$ which contradicts $[U, U]\neq 0.$ Consequently, there exists
$u\in U$ such that $ub\neq 0.$ On the other hand, by Lemma \ref{l4},
there exists a nonzero $\sigma$-ideal $M$ such that $[M,R]\subset U$
and $[M,R]\not\subset Z(R).$ Let $u\in U,\; m\in M$ and $r\in R$; in
view of $[mau,r]\in [M,R]\subset U$ it follows that
\begin{eqnarray*}
  0 &=& a[mau,r]b\\
   &=& a[ma,r]ub+ama[u,r]b\\
  &=& a[ma,r]ub\\
   &=&amarub-armaub\\
   &=&amarub
\end{eqnarray*}
accordingly $$amaRub=0,\;\; \mbox{for all }\;\; m\in M. $$
Similarly, from  $0=\sigma(a)[m\sigma(a)u,r]b$ we find that
$\sigma(a)m\sigma(a)Rub=0$ for all $m\in M$. As $\sigma(M)=M$, then
$ amaRub=0=\sigma(ama)Rub$. Since $R$ is $\sigma$-prime and $ub\neq
0,$ then $ama=0$ and therefore $aMa=0.$ On the other hand
\begin{eqnarray*}
  0 &=& a[m\sigma(a)u,r]b\\
   &=& a[m\sigma(a),r]ub+am\sigma(a)[u,r]b\\
  &=& a[m\sigma(a),r]ub\\
   &=&am\sigma(a)rub-arm\sigma(a)ub\\
   &=&am\sigma(a)rub.
\end{eqnarray*}
Hence $am\sigma(a)Rub=0$ for all $m\in M.$ Since $M$ is a
$\sigma$-ideal we conclude that
$$am\sigma(a)Rub=0 =\sigma(am\sigma(a))Rub.$$
Since $ub\neq 0$, the $\sigma$-primeness of $R$ yields
$am\sigma(a)=0$ thereby $aM\sigma(a)=0.$ Consequently,
$\;aMa=0=aM\sigma(a)\,$ and $\,a=0\,$ by Lemma \ref{l1}. \qed
\vspace*{-0,3cm}
\bl\label{l6}

Let $U$ be a $\sigma$-square closed Lie ideal of a $2$-torsion free
$\sigma$-prime ring $R$ having a nontrivial automorphism $T$
centralizing on $U$ and commuting with $\sigma$ on $U.$ If $x$ in
$U\cap Sa_{\sigma}(R)$ is such that $T(x)\neq x$, then $x$ in
$Z(R).$

\el\vspace*{-0,5cm} \p If $U\subset Z(R)$, then $x\in Z(R).$ Hence
we can suppose that $U\not\subset Z(R)$  and thus $[U, U]\neq 0$
from Lemma \ref{l3}.  Linearizing $[u, T(u)]\in Z(R)$ we obtain $[u,
T(y)]+ [y, T(u)]\in Z(R)$ for all $u,y\in  U.$ In particular, we
have\\$[u, T(u^2)] + [u^2, T(u)]\in  Z(R)$ so that $(u+T(u))[u,
T(u)] \in Z(R).$ Therefore,  $0=[u, (u+T(u))[u, T(u)]]=[u, T(u)]^2.$
Since $[u, T(u)]\in Z(R),$ this forces $[u, T(u)]=0$ for all $u\in
U$. Linearizing this equality, we obtain
\begin{equation}
 [T(u),y]=[u,T(y)]\;\;\mbox{ for all}\;\;u, y\in U.\vspace*{-0,2cm}
\end{equation}
Let $x$ in $U\cap Sa_{\sigma}(R)$ with $T(x)\neq x.$ Replace $y$ by
$2xy$ and replace $u$ by $x$ in $(1),$ using  char$(R)\neq 2,$ we
obtain $0=(x-T(x))[T(x), y]$ for all $y\in U.$ Write $2uy$ instead
of $y$, we get  $(x-T(x))u[T(x), y]=0$ so that $(x-T(x))U[T(x),
y]=0$ and therefore \vspace*{-0,3cm}$$(x-T(x))U[T(x),
y]=(x-T(x))U\sigma([T(x), y])=0\;\;\mbox{ for all}\;\;y\in
U.\vspace*{-0,1cm}$$ Since $T(x)\neq x,$ Lemma \ref{l5} yields that
$[T(x), y]=0$ for all $y$ in $U.$ Hence $[T(x), ry]=[T(x), yr]$ and
thus $[T(x),r]y=y[T(x), r]$ for all $r\in R.$ Replacing $r$ by $ru$,
where $u\in U,$ we find that
$$[T(x), r]uy=y[T(x),r]u=[T(x),r]yu$$ in such a way that $[T(x),
r][u,y]=0$ which proves that $[T(x), r]R[U,U]=0.$ As $U$ is
invariant under $\sigma,$ then $$[T(x), r]R[U,U]=0=[T(x),
r]R\sigma([U,U]).$$ Since $[U, U]\neq 0,$ the $\sigma$-primeness of
$R$ yields $[T(x), r]=0$ so that $T(x)\in Z(R).$ Since $T$ is an
automorphism, it then follows that $x\in Z(R).$ \qed
\\
\\
Now we are ready to state the main theorem of this section.
\vspace*{-0,4cm}
\bt\label{t1}

Let $R$ be a 2-torsion free $\sigma$-prime ring having an automorphism\\
$T\neq 1$ centralizing on a $\sigma$-square closed Lie ideal $U.$
If $T$ commutes with $\sigma$ on $U,$ then $U$ is contained in
$Z(R).$

\et
\p If $[U, U]=0,$ then $U\subset Z(R)$ by Lemma
\ref{l3}. Hence, suppose that $[U, U]\neq 0.$ If $T$ is the identity
on $U,$ then
\begin{equation}
 T([r,x])=[r,x]=[T(r),x],\;
\mbox{ for all}\; r\in R,\,x\in U.
\end{equation}
Let us write $ru$ instead of $r$ in (2), where $u\in U,$ thereby
obtaining
\begin{equation}
T(r)[u,x]=r[u,x]
\end{equation}
Now for any $r'\in R$, If we replace $r$ by $rr'$ in (3) and then
employ $(3),$ we have
$$(T(r)-r)r'[u,x]=0\vspace*{-0,3cm}$$ in such a way that
$\;(T(r)-r)R[u,x]=0 \;\;\mbox{ for all}\;\; u,x\in U, r\in R.$
Since $U$ is invariant under $\sigma,$ we then conclude that
$$(T(r)-r)R[u,x]=0=(T(r)-r)R\sigma([u,x]).$$ From $[U, U]\neq 0,$ it follows that
$T(r)=r$ for all $r\in R$ which is impossible. Hence $T$ is
nontrivial on $U.$ Since $R$ is $2$-torsion free,  then  $T$ is also
nontrivial on $U\cap Sa_{\sigma}(R).$ Therefore, there must be an
element $x$ in $U\cap Sa_{\sigma}(R)$ such
 that $x \neq T(x)$ and $x$ is then in $Z(R)$ by Lemma \ref{l6}. Let $0\neq y$ be in
 $U\cap Sa_{\sigma}(R)$ and not be in $Z(R).$ Once again using
Lemma \ref{l6}, we obtain $T(y)=y.$ But then
\vspace*{-0,3cm}$$T(xy)=T(x)y=xy\;\;\mbox{ so
that}\;\;(T(x)-x)y=0.\vspace*{-0,4cm}$$ Since $x\in\,Z(R),$ it
follows that $\;(T(x)-x)Ry= (T(x)-x)R\sigma(y)=0$  proving $y=0,$ a
contradiction. Hence for all $y$ in $U\cap Sa_{\sigma}(R)$, $y$ must
be in $Z(R).$\\Now let $x$ in $U,$ the fact that $x-\sigma(x)$ and
$x + \sigma(x)$ are in $U\cap Sa_{\sigma}(R)$ assures that both
$x-\sigma(x)$ and $x + \sigma(x)$ are in $Z(R)$ and thus $2x$ in
$Z(R).$ Consequently, $x$ in $Z(R)$ which proves $U \subset Z(R).$
\qed
\vspace*{-0,4cm}
\section{ {\small {\bf  Jordan left derivations on $\sigma$-square closed Lie ideals}} }
\vspace*{-0,4cm}In this section $R$ will always denote a $2$-torsion
free $\sigma$-prime ring and $U$ a $\sigma$-square closed Lie ideal
of $R.$ To introduce the main result of this section, we first state
a few known results which will be used in subsequent discussion.
 \vspace*{-0,4cm}
 \bl\label{l7}

{\em  (\cite{Ash}, Lemma 2.2)} Let $R$ be a $2$-torsion free ring
and let
 $U$ be a square closed Lie ideal of $R.$ If $d\,: \,R\rightarrow
 R$ is an additive mapping satisfying $d(u^2)=2ud(u)$ for all
 $u\in U,$ then \\
\em{(i)} $d(uv + vu)=2ud(v)+2vd(u),$ for all $u,v\in U.$\\
(ii) $d(uvu)=u^2d(v)+3uvd(u) -vud(u),$ for all $u,v\in U.$\\
(iii) $d(uvw + wvu)=(uw + wu)d(v)+3uvd(w)+3wvd(u)-vud(w)-vwd(u),$ for all $u,v, w\in U.$\\
(iv) $[u, v]ud(u)=u[u, v]d(u),$ for all $u,v\in U.$\\
(v) $[u, v](d(uv)-ud(v)-vd(u))=0,$ for all $u, v\in U.$

 \el
 \bl\label{l8}

 {\em (\cite{Ash}, Lemma 2.3)} Let $R$ be a $2$-torsion free ring and let
 $U$ be a square closed Lie ideal of $R.$ If $d\,: \,R\rightarrow
 R$ is an additive mapping satisfying $d(u^2)=2ud(u)$ for all
 $u\in U,$ then \\
\em{(i)} $[u,v]d([u,v])=0,$ for all $u,v\in U.$\\
(ii) $(u^2 v -2uvu + vu^2)d(v)=0,$ for all $u,v\in U.$

 \el
 \vspace*{-0,3cm}
The fundamental result of this section states as follows:
 \vspace*{-0,3cm}
\bt\label{t2}

Let $R$ be a 2-torsion free $\sigma$-prime ring and let $U$ be a
$\sigma$-square closed Lie ideal of $R$. If $d: R\longrightarrow
R$ is an additive mapping which satisfies $d(u^2)=2ud(u)$ for all
$u\in U$, then $d(uv)=ud(v)+vd(u)$ for all $u,v\in U$.

\et
\vspace*{-0,5cm}
\p
 Suppose $[U,U]=0$ and let $u,v\in U$. From $\;d(u+v)^2=2(u+v)d(u+v),$ it follows that $\;2d(uv)=
2ud(u)+2vd(v)-d(u^2)-d(v^2)+2vd(u)+2ud(v),$ in such a way that
$\;2d(uv)=2(ud(v)+vd(u))$ for all $\;u,v\in U.$ As $char R\neq 2,$
then $d(uv)=ud(v)+vd(u).$ Hence we shall assume that $[U,U]\neq 0.$
According to Lemma \ref{l7} (iv) we have \vspace*{-0,2cm}
\begin{equation}
(u^2v-2uvu+vu^2)d(u)=0, \;\; \mbox{for all}\;\; u,v\in U.
\vspace*{-0,2cm}
\end{equation}
Replacing $u$ by $[u,w]$ in (4), where $w\in U,$ and using Lemma
\ref{l8} (i) we obtain $[u,w]^2vd([u,w])=0$ which implies that
$\;[u,w]^2Ud([u,w])=0\;$ for all $\;u, w\in U.$ Let $x, y\in
Sa_\sigma(R)\cap U,$ we have
$[x,y]^2Ud([x,y])=0=\sigma([x,y]^2)Ud([x,y])$ and by virtue of Lemma
\ref{l5} either $\;[x,y]^2=0$ or $d([x,y])=0.$ If $d([x,y])=0,$
applying Lemma \ref{l7} (i) together with char$(R)\neq 2,$ we find
that $d(xy)=xd(y)+ y d(x).$ Now suppose that $[x,y]^2=0$; from Lemma
\ref{l8} (ii) it follows that
 $$(u^2 v -2uvu + vu^2)d(v)=0,\;\;\mbox{ for all}\;\;u, v\in
 U.$$ Linearizing this relation in $u$, we obtain
 $$(uwv + wuv-2uvw - 2wvu + vuw + vwu)d(v)=0,\;\;\mbox{ for all}\;\;u, v, w\in
 U.$$ Replacing $v$ by $[x, y]$ and using Lemma \ref{l8} (i), we
 conclude that
 \begin{equation}
 (-2u[x, y]w - 2w[x, y]u + [x, y]uw + [x,
 y]wu)d([x, y])=0.
 \end{equation}
Write $u[x, y]$ instead of $u$ in (5), since $[x, y]^2=0,$ Lemma
\ref{l8} (i) leads us to
$$
[x, y]u[x, y]wd([x, y])=0,\;\;\mbox{ for all}\;\;u, w\in
 U.
$$ Accordingly, $\;[x, y]u[x, y]Ud([x, y])=0$ for all $u\in
U.$ As $[x,y]\in U\cap Sa_\sigma(R),$ the fact that $\sigma(U)=U$
yields
$$[x,y]u[x,y]Ud([x,y])=0=\sigma([x,y]u[x,y])Ud([x,y])$$
and using Lemma \ref{l5}, either $\;d([x,y])=0\;$ or
$\;[x,y]u[x,y]=0\;$ for all $\;u\in U.$\\If $d([x,y])=0$, then
$d(xy)=xd(y)+ y d(x)$ by Lemma \ref{l7} (i). If $[x,y]u[x,y]=0$ for
all $u\in U,$ then $ [x,y]U[x,y]=0=\sigma([x,y])U[x,y].$ Once again
using Lemma \ref{l5}, we get $[x,y]=0$ and Lemma \ref{l7} (i) forces
$d(xy)=xd(y)+ yd(x).$ Consequently, in both the cases we find that
\begin{equation}
d(xy)=xd(y)+ yd(x), \;\; \mbox{for all}\;\; x,y\in U\cap
Sa_\sigma(R).
\end{equation}
Now, let $u,v\in U$; if we set\\\\

$ \left\{
\begin{tabular}{lll}
$u_1= u + \sigma(u),\;\;\;u_2=u - \sigma(u)$ \\
$v_1= v + \sigma(v),\;\;\;\;v_2=v - \sigma(v)$
\end{tabular}
\right.$  \\\\then we have $2u=u_1+u_2$ and $2v=v_1+v_2 .$ Since
$u_1,u_2,v_1,v_2\in U\cap Sa_\sigma(R)$, application of $(6)$ yields
\begin{eqnarray*}
d(2u2v)&=&d(u_1v_1+u_1v_2+u_2v_1+u_2v_2)\\
 &=&u_1d(v_1)+v_1d(u_1)+u_1d(v_2)+v_2d(u_1)\\
 &+&u_2d(v_1)+v_1d(u_2)+u_2d(v_2)+v_2d(u_2)\\
 &=& 2ud(2v)+2vd(2u).
 \end{eqnarray*}
 As char$(R)\neq 2,$ it then follows that $d(uv)=ud(v)+vd(u)$, for all $u,v \in U.$
 \qed
 \vspace*{-0,2cm}
\bc

 Let $R$ be a $2$-torsion free $\sigma$-prime ring. Then every Jordan left derivation on $R$ is a left derivation on
$R$.

\ec \vspace*{-0,4cm}

\end{document}